\documentstyle[12pt]{article}
\def\DJ{\hbox{D\kern-.8em\raise.15ex\hbox{--}\kern.35em}}

\def\halm{\vrule height1ex width.9ex depth.1ex} 

\def\GL{{\rm GL}}

\def\PGL{{\rm PGL}}

\def\Spin{{\rm Spin}}

\def\SL{{\rm SL}}

\def\SO{{\rm SO}}

\def\SU{{\rm SU}}

\def\db #1 {{\bf#1}}

\def\db#1 {{\bf#1}}
\def\A{{\rm A}}
\def\B{{\rm B}}
\def\C{{\rm C}}
\def\D{{\rm D}}
\def\E{{\rm E}}
\def\F{{\rm F}}
\def\G{{\rm G}}
\def\R{{\rm R}}
\begin{document}
\title{ Rationality
of Almost Simple Algebraic Groups\thanks{A.M.S. 1991
Mathematics Subject 
Classification : Primary 20G15  , Secondary 14G05, 14G10 }}
\author{Nguy\^e\~n Qu\^o\'{c} Th\v{a}\'ng}
\date{}
\maketitle
\begin{center}Department of Mathematics and Statistics, McMaster University, 
Hamilton\\ Ontario, Canada L8S 4K1\\
E-mail : nguyen@icarus.math.mcmaster.ca
\end{center}
\begin{abstract} We prove the stable rationality of almost 
simple algebraic groups, the connected components of the Dynkin diagram 
of anisotropic kernel of which contain at most two vertices.
The (stable) rationality of many isotropic almost simple
groups with small anisotropic kernel and 
some related results over $p$-adic and arbitrary fields are discussed.
\end{abstract}
{\bf Introduction}
\\
Let $G$ be a linear algebraic group defined over a field $k$.
The well-known results of Chevalley and Cartier
showed that if $k$ is an algebraically       
closed field and $G$ is connected and reductive then  $G$ is rational
over $k$ as $k$-variety, i.e., the field $k(G)$ of rational
functions defined over $k$ of $G$ is a pure transcendental  extension
of $k$. However  this is no longer true if $k$ is not algebraically closed
and one of basic geometric problems of algebraic groups over 
non-algebraically closed fields is the rationality problem.
A milder notion of stable rationality (and unirationality) is in sequence :
An irreducible $k$-variety $X$ is $k$-stably rational (resp. $k$-unirational) 
if there is an affine $k$-space $A$ such that $X \times A$ is 
$k$-birationally equivalent to an affine $k$-space (resp. such that there is
a surjective $k$-morphism $A \to X$). In general, it is difficult 
to verify
if a given $k$-group (or $k$-variety) is rational (or irrational).
We refer the readers to [Ch], [CT], [MT], [M1-2], [P], [V], [VK]
and references thereof for various
problems and progress related with the rationality problem.
\\
Up to now there is no general critetion to decide which almost simple groups
are stably rational over the field of definition by looking at their
Dynkin diagram, except the trivial cases of split and quasi-split groups.
Quite recently it became known that for many division algebras
of degree 4 (or greater) over $k$, there are many examples of $k$-groups,
isotropic or not, related 
with them which are
not stably rational over $k$. The purpose of this note is to show 
that this result can be used to get such a general criterion.
In particular we show that many almost simple groups with
relatively big $k$-rank  and the degree of the related division
algebra is $\le 3$ are $k$-stably rational. Hence in certain sense, our
results are  optimal. More precisely the main result of the paper 
is the following.  
\\
\\
{\bf Theorem.} {\it Let G be an  
almost simple algebraic group over a field k. Let $m(G)$ be the
maximal number
of vertices of connected components of the Dynkin diagram
of anisotropic kernel of G.}\\
$a)$ {\it If $m(G) \le 2$, then G is rational or stably rational over k.}\\
$b)$ {\it The number $2$ in $a)$ is best possible.  For any natural number n, 
there exist non stably rational groups G and fields k with  
$m(G)=4n-1$.} 
\\
\\
In particular over certain
"nice" fields, such as local ($p$-adic or real) fields, 
many isotropic almost simple
$k$-group are $k$-(stably) rational. 
The method of the proof is based 
on a detailed  analys of the Tits index of the groups under consideration.
\\
{\it Notation.} For an almost simple group $G$ defined over a field $k$
which has characteristic either 0 or relatively big, e.g. relatively
prime with the order of
the center $Cent(G)$ of $G$.
Let $S$ be a maximal $k$-split torus of $G$, $T$ a maximal $k$-torus
of $G$ containing $S$. If dim $T=n$, we denote by 
$$\Delta = \{\alpha_1, \ldots, \alpha_n \}$$ a basis of simple roots
for the root system $\Phi$ of $G$ with respect to $T$.  We may consider
the relative root system $_k \Phi$ of $G$ relatively to $S$ and let
$_k \Delta$ a basis of $_k \Phi$ compatible with $\Delta$.
For $1 \le i \le n$   we denote by $S_i$ the standard
$\bar k$-split torus corresponding to the root $\alpha_i$. We 
denote by $x_\alpha (t)$ the multiplicative one-parameter unipotent subgroup
(resp. $h_\alpha (t)$ the
multiplicative one-parameter diagonal subgroup) of $G$ 
corresponding to a root $\alpha \in \Delta$ where we keep 
the same notation used in [St]. For $\alpha =\alpha_i$ we denote
$x_i(t) = x_{\alpha_i}(t)$, $h_i(t) =h_{\alpha_i}(t)$, $1 \le i \le n$
and $X_i$ the image of $x_i$ in $G$.
In particular, if $G$ is simply connected then $T$ is the direct product
of the images of $h_i :=h_{\alpha_i}$, $1 \le i \le n.$
We use intensively the notion and results of Tits' classification theory
of almost simple algebraic groups as presented in [Ti1] and refer
also to [BT] for other notions in algebraic groups.
We often identify a simple root with the vertex 
representing it in the Tits index.
\\
I would like to thank C. Riehm for providing me with the up-dated
version of [M2] and Department of Mathematics and Statistics, McMaster
University for support during the preparation of this paper.
\section {Some general useful facts}
{\bf 1.1.} Let $G$ be a connected reductive $k$-group, $S$ a 
maximal $k$-split torus of $G$.  The Bruhat decomposition for $G$
(see [BT, Section 4])
implies that $$G \simeq  Z_G(S) \times {\bf A}$$ as varieties,
where {\bf A} is  an affine space defined over $k$. Thus the study 
of rationality of $G$ is reduced, in certain sense, to that of $Z_G(S)$.
Namely  $G$ is stably rational 
over $k$ if and only if $Z_G(S)$ is and if
$Z_G(S)$  is $k$-rational, then so is $G$.
\\
However in certain cases the group $Z_G(S)$ is hard to handle with and we are forced
to find a substitute, which can be studied easier. In many cases
it is possible to do so. Namely let
$S_0$ be  a nontrivial $k$-subtorus of $S$. Another version of Bruhat
decomposition says that $$G \simeq Z_G(S_0) \times {\bf A},$$ the direct
product 
of $Z_G(S_0)$ with an affine space {\bf A} over $k$. 
\\
 Therefore we are reduced to 
studying the connected reductive $k$-groups $Z_G(S_0)$. The problem here
is to choose a "nice" torus $S_0$ so that we can prove the
rationality or stable rationality
of $Z_G(S_0)$, which is possible if the $k$-rank is relatively
big. First we need the following simple
but very useful  observation.\\
\\
{\bf 1.2. Proposition.} [DT] {\it Let $S_{\theta}$ be a standard k-split torus
of $G$ and $Z_G(S_{\theta}) =S_{\theta}T_0 H$ (almost direct product),
where $T_0$ is a k-torus, H a semisimple k-subgroup of G. Then the Tits
index of H is obtained from that of G by remowving all vertices not
belonging to the preimage $\tilde \theta$ of $\theta \cup \{0\}$
under the restriction map $\Delta \to~ _k\Delta \cup \{0\}.$ Moreover
$T_0$ is anisotropic and $ST_0 = (Z_T(H))^0$.}  
\\
\\
{\bf Remark.} 
The equality in the last statement  is not in [DT] but it is clear by 
comparing the dimension of both side and by making use of the previous part
of the proposition.  In particular it shows that if $Z_T(H)$ is connected
then $ST_0$ contains the center of $H$.
\\
Another interesting remarks are the following observations due to Tits.
(See [Ti2] and [Se].)
In fact, $S$ can be replaced by a standard $k$-split torus,
but we do not need this here.\\
\\
{\bf 1.3. Proposition.} {\it Let G be 
an almost simple adjoint $k$-group with
a maximal k-split torus S. Then the center of $Z_G(S)$ is connected
if G is either adjoint or simply connected.}
\\
\\
{\bf 1.4. Proposition.} {\it Let  G be an almost simple k-group,
$S_{\theta}$ a standard $k$-split torus of G. Let }
$$Z_G(S_{\theta}) = S_{\theta}T_0 H.$$
{\it Then $S_{\theta}T_0$  is a quasi-split k-torus
hence also cohomologically  trivial.}\\
\\
The following result essentially is due to Ono - Rosenlicht (see [O]).\\
\\
{\bf 1.5. Proposition.} {\it If T' is a central k-torus of a connected
reductive k-group G, which is quasi-split torus over k, then there is
a rational k-cross section $G/T' \to G$. In particular the k-variety 
G is birationally equivalent to the product $T' \times (G/T')$.}
\\
\\
>From above we see that it is essential to know the group
$Z_G(S)/(ST_0)$ (which we call the {\it semisimple
anisotropic quotient} of $G$) if we want to know the rationality
property of $G$. In the next section we examine various computations of
this group. The following remark is useful in the sequel.
\\
\\
{\bf 1.6. Remark.} {\it If $\pi : G \to G'$ is a central k-isogeny  
and $S' = \pi(S)$ the image of maximal k-split torus of G, then $\pi$
induces a central isogeny}
$$ Z_G(S)/S \to Z_{G'}(S')/S'$$ 
{\it of semisimple anisotropic quotients.}
\\
\\
\section {Some computations related with centralizers of split tori}
{\bf 2.1.}
We keep the above notation and we assume that
 $G$ is an almost simple $k$-group. Let $P_i$,
$1 \le i \le s$, be  $k$-groups with $F_i$
 a central
$k$-subgroup of $P_i$. Assume that all $F_i$ are $k$-isomorphic. By a suitable
factoring out a central $k$-subgroup of the
 direct product of $P_i$ we will obtain
an almost direct product $P_1'... P_s'$ with the property that the
set-theoretic intersections $P'_i \cap P'_j$ are all equal and $k$-isomorphic
to $F_i$. We call such a group the {\it product of $P_i$ with glued central
 subgroups $F_i$.}
\\
\\
{\bf 2.2. Proposition.}
{\it Let G be an almost simple k-group of type $\A_n$ with
k-rank $r>0$.}\\
{\it a) If G is an inner form, the anisotropic semisimple 
quotient $Z_G(S)/S$ is the product
of $k$-conjugate 
almost simple anisotropic k-groups of type $^1\A_{d-1}$ with glued center.}
\\
{\it b) If G is an outer form, $Z_G(S)/ST_0$ is the product of  
anisotropic groups
of type
$^1\A_{d-1}$ with an anisotropic k-group of type $^2\A_{n-2rd}$ with glued
central
subgroup of order dividing d.}\\
\\
{\it Proof.}  
$a)$ First we begin with simply connected
isotropic groups of type $^1\A_n$. From [Ti1] 
we know
that such groups have the following Tits index
\\
\\
$$\bullet^1-- \cdots -- \bullet^{d-1} -- \odot^d -- \bullet
--\cdots -- \bullet ^{rd-1} -- \odot^{rd} -- \bullet -- \cdots --\bullet ^n$$
\\
\\
Let $S_i$ be the standard $k$-split torus corresponding to the 
isotropic vertex $i$, $i = d, 2d,\ldots, rd$, where $d $ is the index
of the division $k$-algebra $D$ related with this type
and $r$ denotes the $k$-rank of $G$. Then 
$$S = \prod _i S_i$$ is a maximal $k$-split torus of $G$. 
We have $Z_G(S) = S H$, where $H= \prod H_j$ is a semisimple $k$-group which
is an almost direct product of anisotropic $k$-groups $H_j$ of type 
$^1\A_{d-1}$ (see Proposition 1.2). 
\\
We show that all these groups $H_j$ are in fact
$k$-isomorphic to the simply connected almost simple $k$-group
$G_0$ of type $^1\A_{d-1}$ defined by 
$G_0(k) = \SL_1(D)$. Indeed, we may assume that in certain basis, 
the maximal $k$-split torus $S(k)$ consists of
all diagonal matrices from $\GL_n(D)$ with coefficients from $k^*$
and of determinant 1 :
$$S(k) = \{diag(t_1,\ldots,t_n) :~ t_1\cdots t_n =1,
\in k^*, 1 \le i \le n \},$$
thus 
$$Z_G(S)(k) = S(k)\{diag(d_1,\ldots,d_n) ~: d_i \in \SL_1(D), 
1 \le i \le n \}.$$ 
It follows that all the groups $H_j$ above are $k$-conjugate (i.e. conjugate
by elements from $G(k)$). Therefore
if $G'$ is a quotient of $G$ by a central $k$-subgroup
then we can form the corresponding centralizer of a maximal $k$-split torus
$S'$
of $G'$, which is an almost direct product of $S'$ and $k$-conjugate
almost simple $k$-groups $H'_j$ of type $^1\A_{d-1}$;
and they are the homomorphic image of the simply connected
almost simple $k$-groups $\tilde H_j$ of type $^1\A_{d-1}$ such that
$\tilde H_j (k) \simeq \SL_1(D)$. Some tedious computations 
shows that $S$ contains
the products $z_i z_j$ of two generators $z_i$ and $z_j$
of centers of the groups $H_i$
and $H_j$, respectively.
>From this and 1.6, the assertion  $a)$ follows.\\
$b)$ Let $l$ be the separable quadratic extension of $k$ over which 
$G$ becomes of inner type. Assume first that $G$ is simply connected.
>From [Ti1] we know that 
the Tits index of $G$ is as follows 
\newpage
\noindent
\hspace*{0.5cm}$\bullet \cdots \bullet -- \odot--\bullet \cdots \bullet -- \odot
--\bullet \cdots \cdots \bullet --\odot--\bullet \cdots \bullet \cdots$\\
\\
\hspace*{2.6cm}$\updownarrow \hspace{3.1cm} \updownarrow \hspace {3.59cm}\updownarrow \hspace{3cm}\Bigg )$
\\
\\
\hspace*{0.5cm}$\bullet \cdots \bullet -- \odot--\bullet \cdots \bullet -- \odot
--\bullet \cdots \cdots \bullet --\odot--\bullet \cdots \bullet \cdots$\\
\\
\\
Denote by $G(\Psi)$ the semisimple regular subgroup of $G$ generated
by the root subgroups $X_i$, $\alpha_i \in \Psi$, where $\Psi$ is
a subset of $\Delta$.
By Proposition 1.2 in order  
to compute the intersection  $ST_0 \cap H$ we are reduced to computing
the intersection of the torus 
$$S_dS_{n-d+1} S_{2d}S_{n-2d+1} \cdots S_{rd}S_{n-rd+1}$$ 
with the
semisimple subgroup $$G'_1\cdots G'_r A,$$
where $G'_i$, $ 1 \le i \le r$ is the semisimple $k$-group
$G(\Psi_i)$ with the root system generated by the basis
$$\Psi_i = \{
\alpha_{(i-1)d+1},\ldots,\alpha_{id-1},
\alpha_{n-(i-1)d},\ldots,
\alpha_{n-id+2}\},$$
and $A$ is the group $G(\alpha_{rd+1},\ldots,\alpha_{n-rd})$.
Then as in the part $a)$ the part $b)$ follows. The general case also follows
from this as in $a)$ by making use of 1.6.
\halm
\\
\\
The case $d \le 3$ is of special interest to us
due to the following results.\\
\\
{\bf 2.3. Proposition.} {\it $a)$ Let G be an almost simple k-group
of type $^2\A^{(1)}_{n,r}$ (i.e. d=1). With above notation, 
$ST_0$ contains the center of H.}\\ 
{\it $b)$ If G is above and n is even then G is k-rational.}\\
{\it $c)$ 
If G is as above and 
n is odd, then any  almost  simple k-groups which is isogeneous
to G is also k-birationally isomorphic (as varieties) to G.}\\
\\
{\it Proof.} $a)$ We assume first that 
$G$ is simply connected. The general case follows from this since
if $\pi : G \to G'$ is a central $k$-isogeny, then $S' :=\pi (S)$ is
a maximal $k$-split torus of $G'$ and $\pi (Z_G(S)) = Z_{G'}(S')$
and we use 1.6.
For simplicity we assume that $r=1$ and we give a 
complete computation in this case. From Proposition 1.2 it follows that
we have only to check that 
$$Cent(H) \subset S_1S_d,\eqno{(1)}$$
where
$H =G(\alpha_2,\ldots,\alpha_{n-1}).$ 
>From above we see that for an element
$t \in T$,
$$t = \prod_{1 \le i \le n} h_i(t_i)$$
is in $Z_T(H)$ if and only if $t$ commutes with all one-parameter
unipotent subgroups $X_i$ for $2\le i \le n-1$. 
Hence we have the following system of equations for $t_i$ :
\\                                                
\\
\[ \left\{  \begin{array}{l}
t_2^2 = t_1t_3,
\\
t_3^2=t_2t_4,
\\
\vdots
\\
t_{n-1}^2 =t_{n-2}t_n.
\end{array}
\right. \]
\\
\\
One checks that $$t_i= t_2^{i-1}/t_1^{i-2}, 3 \le i \le n,$$ 
while the center
of $H$ is generated by
$$h_2(\zeta)\cdots h_{n-1}(\zeta^{n-2}),$$ where $\zeta$ is a primitive
$(n-1)$-root of unity. Hence  (1) is verified and $a)$ follows.
\\
$b)$ Note that (with notation as above), the isogeny $\pi$ 
induces an isogeny (denoted by the same symbol)
$$\pi: Z_G(S)/ST_0 \to Z_{G'}(S'T'_0)/S'T'_0$$
between the anisotropic semisimple quotients.
Now by $a)$ $ST_0$ contains the center of $H$, the corresponding anisotropic 
semisimple quotient is an {\it adjoint} group, hence the above 
induced isogeny is in fact a $k$-isomorphism. 
\\
Now $b)$
follows from $a)$, results of Section 1 and a result of [VK]
that  any adjoint $k$-group of type $\A_{2m}$, $m \ge 1$, is $k$-rational.
\\
From Section 1 it follows that 
$$Z_G(S) \simeq ST_0 \times (Z_G(S)/(ST_0),$$
and from  $ST_0 \simeq S'T'_0$ it follows that $G$ and $G'$ are
birationally isomorphic over $k$ hence $c)$.
\halm
\\
\\
From 2.2 we derive the following result regarding $d \le 3$.\\
\\
{\bf 2.4. Proposition.} {\it Let k be a 
field with a division algebra D of degree $d \le 3$, $G_0$ an 
anisotropic semisimple
k-group which is an almost direct product of k-groups of type $\A_{d-1}$ isogeneous
over k to either $\SL_{1,D}$ or $\PGL_{1,D}$ such that any simply
connected factor contains the center of the other.
Then $G$ is stably rational over k. In particular, if k has a unique up to
isomorphism quaternion division algebra (e.g. k is a local field),
then any such almost direct 
product of anisotropic k-groups of type $\A_1$ is stably rational over k.}\\ 
\\
{\it Proof.}
Since $\PGL_{1,D}$ has no center and is $k$-rational, 
we may assume that all almost simple factors of
$G$ are isomorphic to $\SL_{1,D}$ and they have common center (i.e.
product of groups of type $\A_{d-1}$ with glued center). Let the
number of almost
simple components of $G$ 
be $r$. Then from 2.2 we see that for the group $G_1$
with $G_1(k) = \SL_{r+1}(D)$
and $S$ a maximal $k$-split torus of $G_1$ we have
$$Z_{G_1}(S) = S G_2,$$
and $$ Z_{G_1}(S) /S \simeq  G.$$
Since the group $\SL_{1,D}$ is $k$-rational by assumption on $d$,
the group $\SL_{n,D}$ is
also for any $n$ (see [M1], [V]), and it follows that $G$ is stably
rational over $k$.
\halm
\\ 
\\
{\bf 2.5. Remarks.} $a)$ Another close
formulation of the proposition is as follows.
\\
\\
{\it Let G be an almost simple k-group with the following Tits index}\\
\\
$$ \bullet -- \odot -- \bullet -- \odot \cdots \cdots \bullet -- \odot -- \bullet,$$
\\
{\it i.e. of type $^1\A^{(2)}_{2r+1,r}$, or the following}
\\
$$\bullet -- \bullet -- \odot -- \bullet -- \bullet \cdots \bullet -- \bullet
-- \odot -- \bullet -- \bullet,$$
\\
{\it i.e. of type $^1\A^{(3)}_{3r+2,r}$.
Then G is stably rational over k, and it is
rational over k if it
is adjoint.}
\\
\\
{\it Proof.} Only the case of adjoint groups needs a proof. But in this case
we make use of Proposition 1.3 and notice that adjoint groups of type
$\A_s$, $s \le 2$ are rational.\\
\\
$b)$ It is not clear if $G$ (above)  is always rational, though 
in general the cancelation of rational varieties does not hold.
\\
\\
$c)$ The well-known examples of non-rational 
almost simple groups are related with certain division algebras 
of index 4 or greater. The results of Merkurjev [M1-2] show that if $G$ is
simply connected of type $^1\A_n$ such that $G(k) = \SL_n(D)$ where
the index of $D$ is divisible by 4 then $G$ is not stably rational over $k$
and also there exist adjoint groups of type $\A_3$ which are not 
stably rational.
It is not known if there are non-rational 
(or non stably rational) 
adjoint groups of type $\A^{(d)}_{n,r}$ related with division algebra 
$D$ of degress $d\le 2$.
\\
\section {Rationality of almost simple groups over local fields}
{\bf 3.1.} First we asssume that $k$ is a $p$-adic field. We have the 
following result about rationality of almost simple $k$-groups over $k$.\\
\\
{\bf 3.2. Theorem.} {\it Assume that G is an 
almost simple group over a p-adic field k.}\\ 
{\it $a)$ If G is adjoint then G is k-rational.}\\ 
{\it $b)$ If G is of 
type different from $ ^1\A^{(d)}$ $(d \ge 4)$, simply connected types
$^1\D^{(2)}_{2r+3,r}$ or $^2\D^{(2)}_{2r+2,r}$ (r odd) 
then G is rational or stably rational over k.}\\
\\
{\it Proof.} $a)$ From Section 1 we know that if $S$ is a maximal $k$-split
torus of $G$ then the center of $Z_G(S)$ is connected. A theorem of
Kneser says that anisotropic semisimple $k$-group is necessary an almost
product of groups of type $^1\A$. Since the groups $\PGL_{1,D}$ are
$k$-rational, $a)$ follows.
\\
$b)$ First we assume that $G$ is of classical type. If $G$ is of type A 
then from [Ti1] and the assumption we deduce that it is of outer type
$^2\A_{2r,r}$ or $^2\A_{2r+1,r}$. From Section 1 it follows easily that
$G$ is rational. If $G$ is simply connected of type B or $^1\D^{(1)}$, 
then the assertion follows from the fact that the Spin group
of an isotropic non-degenerate quadratic form is rational (see [P]). If
$G(k)$ is given by the special orthogonal or unitary group of some
quadratic or skew-hermitian form, $G$ is rational due to Cayley 
transformation. 
\\
Now we assume that $G$ is simply connected of type $\D^{(2)}$.
>From [Ti1] it follows that $G$ is of type $^1\D_{2r,r}$
with the following Tits index
\\
\\
\hspace*{9.7cm}$\bullet^{2r-1}$\\
\hspace*{2cm} $\bullet^1--\odot^2--\cdots--\bullet^{2r-3}--^{2r-2}\odot 
$\\
\hspace*{9.7cm}$\odot^{2r}$\\
\\
or $^2\D_{2r+1,r}$ with the Tits index
\\
\\
\hspace*{9.7cm}$\circ^{2r}$\\
\hspace*{2cm} $\bullet^1--\odot^2--\cdots--\odot^{2r-2}--^{2r-1}\bullet
\hspace{0.45cm}\updownarrow$\\
\hspace*{9.7cm}$\circ^{2r+1}$\\
\\
In either case, for a maximal split $k$-torus $S$ of $G$
we have
$$Z_G(S) = ST_0 H,$$
where $H$ is the direct product of simply connected 
anisotropic groups of type $\A_1$. One can check that the anisotropic
quotient $Z_G(S)/ST_0$ is an almost direct product of groups
isomorphic to $\SL_{1,D}$ with common centers, or just direct product
of adjoint groups of type $\A_1$. The latter groups are known
from above (see 2.3) to be stably rational over $k$.  Thus $G$ is also.
Moreover, for any central $k$-isogeny $\pi : G \to G'$, the anisotropic 
semisimple quotient of $G'$ is the image of that of $G$,
hence is also the product of groups of type $\A_1$
with common center hence is stably rational as above. This covers also
the case
of almost  simple groups of type $\D_{2m}$ 
with center of order 2
which are obtained by factoring the Spin group $\tilde G = \Spin(\Phi)$ 
by the subgroups
$\{1,z\}$, $\{1,z'\}$, where the center of $\tilde G$ is
$\{1,z,z',zz'\}$ and  $\SO(\Phi)$ (or $\SU(\Phi)$) is
$\simeq \Spin(\Phi)/\{1,zz'\}$. (Note that the roles of the two elements
$z$ and $z'$ in certain cases are {\it not} symmetric.)
\\
Now we consider the case $^2\D^{(2)}_{2r+2,r}$, $r$ is even, $r=2s$. 
First we assume $G$
simply connected. The Tits index of $G$ is as follows
\\
\\
\hspace*{9.5cm}$\bullet^{2r+1}$\\
\hspace*{2cm}$\bullet^1 --\odot^2--\cdots--\bullet^{2r-1}--^{2r}\odot
\hspace{0.75cm}\updownarrow$\\
\hspace*{9.5cm}$\bullet^{2r+2}$
\\
\\
Note that $$Z_G(S) = SG_1 G_3\cdots G_{2r-1} A A',$$
where $G_i$, $1 \le i \le 2r-1$ (odd) are anisotropic simply connected
$k$-group of type $\A_1$, $A A' = \R_{l/k}(A'')$ with a quadratic
extension $l$ of $k$, over which $G$ is of inner form, and $A''$ is
an anisotropic $l$-group of type $\A_1$. We may assume that all $G_i$
are identified with the $k$-group $\SL_{1,D}$ for a unique (up to isomorphism)
quaternion division algebra $D$ over $k$. We notice that the
product $H = G_1 G_3 \cdots G_{2r-1} AA'$ is stably equivalent to the group
$\R_{l/k}(H')$, where $H'=G_1G_3\cdots G_{r-1} A$. (Here $G_i$ are
$k$-groups but considered as $l$-groups.) We can give a similar 
interpretation for the center of $H,H'$.  Since the product $P'$ of
$l$-groups $G_i$ ($1 \le i \le r-1$) and $A$ with glued centers
are stably rational over $l$, it follows that $P= \R_{l/k}(P')$
is stably rational over $k$. But one checks that $P \simeq 
Z_G(S)/S $ (the same computations as above). Thus $G$ is also 
stably rational over $k$. The same is also true for the factor
groups $\Spin(\Phi)/\{1,z\}$, $\Spin(\Phi)/\{1,z'\}$ (see notation
above).
\\
Finally we consider the case of exceptional groups. The non trivial cases
are groups of type $^1\E^{16}_{6,2}$ and $\E^{9}_{7,4}$ and we assume that
they are simply connected.
For
the groups of type $^1\E^{16}_{6,2}$, 
we have the Tits index of $G$ as follows
\noindent
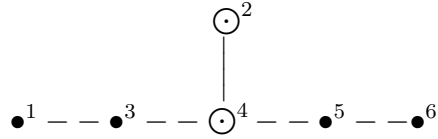
\begin{center}\hspace*{0.09999cm} $\bigodot ^2$\\
$\bigg |$\\
$\bullet^1--\bullet^3--\bigodot^4--\bullet^5--\bullet^6$
\end{center}
and let $S_2$ be the standard split torus corresponding to
the root $\alpha_2$. Then we have 
$$Z_G(S_2)  = S_2 A,$$
where $A$ is a simply connected of type $^1\A^{(3)}_{5,1}$. There is a
division $k$-algebra $D$ such that $A(k) = \SL_{2}(D)$.  Since
the $k$-group $L$ of type A defined by $L(k) = \SL_1(D)$ has rank 2
hence is rational, the group $A$ is also $k$-rational. Now one can check 
that  
$$S_2 \cap A = \{1\},$$ thus from Section 1  we know that $G$ is also
$k$-rational.
\\
The case of groups of type $\E^{9}_{7,4}$ with the Tits index
\begin{center}
\hspace*{0.25cm}$\bullet^2$\\
\hspace*{0.07cm}$\bigg |$\\
$\bullet^7--\odot^6--\bullet^5--\odot^4--\odot^3--\odot^1 \hspace{1cm},$
\end{center}
is reduced to the case of groups
which are product of groups of type $\A_1$ with glued center.
\halm
\\
\\
{\bf 3.3. Remark.} 
The exclusion of groups of type  $^1\A^{(d)}$ with $d \ge 4$ is necessary
(see the end of Section 2). However it is not clear if it is so
regarding the groups of type D in the above proposition.
\\
\\
{\bf 3.4.} 
Now we assume that $k = \db R $. Recently Chernousov [Ch]
has proved that if $G$ is an anisotropic semisimple \db R -group with
no factors of type $\E_6,\E_7, \E_8$ then $G$ is stably rational over \db R .
The main idea there (as in [M2]) is to use the group of similarity factors 
of the
forms involved, which goes back to [T1-3], where we considered the
problem of weak approximation in a close relation with the problem of rationality.
In view of results above, we can generalize the result of Chernousov
as follows.\\
\\
{\bf 3.5. Proposition.} {\it  Let G be a semisimple \db R -group with no 
anisotropic factors of types $\E_i$, i=6,7,8. Then G is stably rational
over \db R .} \halm
\\                
\section {Rationality over arbitrary field 
of some isotropic almost simple groups}
One may notice that some of arguments used above can be applied to a 
more general situation. In fact we have the following \\
\\
{\bf 4.1. Theorem.} {\it Let G be an isotropic  almost simple k-group of
one of the following  types
$^1\A^{(d)}, d \le 3$; $\B$; $\C^{(d)}_{n,r}, n-dr \le 2, d \le 2$;
$\D^{(d)}_{n,r}, d \le 2$, with $n \le r+2~ (d=1)$, $n \le 2r+1~ (d=2)$.
Then G is either rational or stably rational over k.}\\
\\
{\it Proof.} We need only to prove the cases C and D. For the case of type C
we make use of Propositions 1.3, 1.5. For the case of type D, we take the
centralizer of the standard 1-dimensional split torus $S_0$ of $G$ corresponding
to the last circled vertex (in the case $n =2r$) or vertices (in the case
$n=2r+1$) (see the Tits indices drawn above for these groups)
to get the group $$Z_G(S_0) = S_0 T_0 A,$$
where $A$ is a $k$-group of type $^1\A^{(2)}$, which has the factor
group $A/(A \cap S_0T_0)$ of the same type hence  stably
rational as shown above (see Sec. 2). Hence from Section 1 we know
that $G$ is also stably rational. The other cases are treated in a similar
way. \halm
\\
\\
Now we focus our attention  
to exceptional groups and we have the following result.
\\
\\
{\bf 4.2. Theorem.} {\it Let G be an isotropic almost simple group
of exceptional type over a field k.}\\ 
{\it $a)$ If G is of type $\D_4, \F, \G$ then G is k-rational.}\\
{\it $b)$ Let G be one of the following types : $^1\E_6$; $^2\E_{6,r}$
$r \ge 2$ and G is adjoint if G is of type $^2\E^{16"}_{6,2}$;
$\E^{66}_{7,1}$ (simply connected);
$\E_{7,r}, r \ge 2$ (and G is simply connected if
G is of type $\E^{31}_{7,2}$); 
$\E_{8,r}, r \ge 3$. Then G is either k-rational
or stably rational over k.}\\
{\it $c)$ If G is of type $\E^{78}_{7,1}$, then any 
k-group which is k-isogeneous to G is also  birationally isomorphic
to G over k 
as k-varieties.}\\
\\
{\it Proof.} $a)$ The proof follows from results of Section 1.\\
$b)$ For the group of type $^1\E_6$ we consider only the group 
$^1\E^{(28)}_{6,2}$ with Tits index
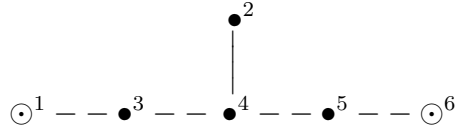
\begin{center}\hspace*{0.08cm} $\bullet^2$\\
$\bigg |$\\
$\odot^1--\bullet^3--\hspace{0.08cm}\bullet^4--\bullet^5--\odot^6$
\end{center}
since for the other, the proof is the same as in $p$-adic
case above. Let $\tilde G$ (resp. $\bar G$) be the simply connected
covering (resp. adjoint group) of $G$ and $\tilde S$ a maximal $k$-split
torus of $\tilde G$. It is well-known that 
$$Z_G(\tilde S) = S D,$$ where $D$ is  the Spin group of a quadratic form $f$
which is a norm form of a division Cayley algebra, i.e., $f$ is a
Pfister form. By a result of [M2] (Prop. 7) the adjoint group of $\SO(f)$
is stably rational over $k$. Hence $\bar G$ is stably rational over $k$.
Now we can check that $\tilde S$ contains the center of $D$. Thus $\tilde G$,
being birationally equivalent to $\bar G$, is also stably rational over $k$.
\\
Let $G$ be simply connected  of type $\E_{6,2}^{16'}$ with the following
Tits index 
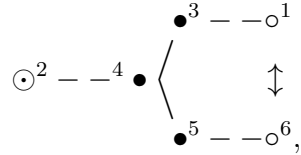
\begin{tabbing}
\hspace*{7cm}
\=$\bullet^3 -- \circ^1$
\\
\hspace*{5cm}$\odot^2 -- ^4\bullet\Bigg \langle\hspace{1.1cm} \updownarrow$\\
\>$\bullet^5 -- \circ^6,$
\end{tabbing}
>From this we see that the anisotropic kernel of $G$ is of type $^2\A^{(1)}_3$,
which is also the anisotropic kernel $A$
for the simply connected semisimple group of 
type $^2\A_5$ with the root system spanned by $\Delta \setminus \{\alpha_2\}$.
By a result of [M2] (Prop. 8), the adjoint group of $A$ is stably rational
and one checks easily that for a maximal $k$-split torus $S$
of $G$, the center of $A$ is contained in $S$. Therefore 
$G$ and its adjoint group are stably rational.
\\
The case of type $^2\E^{}_{6,2}$ is considered similarly as above : If 
$\bar S$ (resp.
$\bar G$) 
has meaning as above, then $\bar S$ contains the center of $Z_{\bar G}(\bar S)$
hence $\bar G$ is rational. 
\\
Now we assume that $G$ is of type $\E_7$. We claim that
if $G$ is simply connected of type $\E^{66}_{7,1}$ then $G$ is rational.
Indeed, for $S$ a maximal $k$-split torus of $G$, we have
$$Z_G(S) = S H,$$
where $H$ is a simply connected $k$-group of type $\D_6$, $H
\simeq \Spin(\Phi)$ for some form $\Phi$. One checks that
we have $$Z_G(S)/S \simeq \SU(\Phi),$$ thus is rational and so is $G$.
The only other non-trivial case
considered here is  the case of type $\E_{7,2}^{31}$ with the Tits index
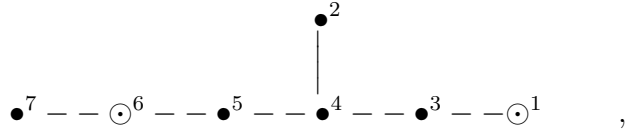
\begin{center}
\hspace*{0.25cm}$\bullet^2$\\
$\bigg |$\\
$\bullet^7--\odot^6--\bullet^5--\bullet^4--\bullet^3--\odot^1 \hspace{1cm},$
\end{center}
Assume that $S$
is a maximal $k$-split torus of the simply connected group $G$   
of this type. We have  
$$Z_G(S_6) = S_6 (\tilde A_1 \times \tilde D_5),$$
where we denote by $\tilde A_1$ (resp. $\tilde D_5$) the group of type
$\A_1$ (resp. $\D_5$) with the root system spanned on
$\alpha_7$ (resp. $\{\alpha_1, \ldots, \alpha_5\}$). We have
$$Z_G(S_6)/S_6 = G_1 G_2,$$
where $G_1$ (resp. $G_2$) is simply connected of type $\A_1$
(resp. $G_2 \simeq \SO(f)$, $f$ is an isotropic non-degenerate
quadratic form in 10 variables over $k$ with Witt index 1). 
Here we have $$G_1 \cap G_2 = \{\pm 1\}.$$ To see this
one just needs to compute the intersection $S_6 \cap (\tilde A_1 
\times \tilde D_5)$ and  the factor $$F/(S_6 \cap F),$$
where $F = Cent (\tilde A_1 \times \tilde D_5)$.
\\ 
Now we want to show that $H =G_1 G_2$ is rational over $k$. Let $S_1$
be a (unique up to conjugacy over $k$) $k$-split torus of
$G_2$.  As above, we show that $Z_H(S_1)$ is rational over $k$. 
We have $Z_H(S_1)= S_1 G_1 G_3,$ where  $G_3 = \SO(f_0)$, where
$f_0$ is the anisotropic part of $f$. We can check that
in fact we have a direct product decomposition
$$Z_H(S_1) = S_1 \times G_1 \times G_3,$$
which is clearly rational over $k$. \\
$c)$ 
It rests to consider the case $\E_{7,1}^{78}$. Let $S$ be a maximal
$k$-split torus of $G$, where we assume that $G$ is simply connected. 
One checks
that $S$ contains the center of the anisotropic kernel of type $\E_6$
of $G$. Thus
$$Z_G(S)/S \simeq Z_{\bar G}(\bar S)/{\bar S}$$ and we are done.
\halm
\\
\\
{\bf 4.3.} We say that a segment in the  Tits index of an almost simple
$k$-group $G$ is $black$ (resp. $white$)
if it consists of only black (resp. white, i.e., distinguished)
vertices. The $length$
of a segment is the number of vertices it contains.
We say that a segment
is {\it defined over} $k$, if the almost simple subgroup of $G$
with root system spanned on 
this segment is defined over $k$. In other words, the black segments
are the connected components of the Dynkin diagram of anisotropic
kernel of $G$. From results proved above we derive
the following main result of this paper. \\
\\
{\bf 4.4. Theorem.} {\it Let G be an almost simple k-group and $m(G)$
be the maximal length of the black segments of its Tits index defined over k.}
\\
$a)$ {\it If $m(G) \le 2$ then
G is either rational or stably rational over k.}\\
$b)$ {\it The number $2$ in $a)$ is best possible. For any natural number
n there exist non stably rational groups G and fields k with 
$m(G)=4n-1$.}\\
\\
{\it Proof.} 
$a)$ It follows from above results and the Tits classification of
indices [Ti1]. \\
\\
$b)$ The number 2 in the above theorem is the best 
possible since [M1] shows that it fails if 2 is replaced by 3.  Namely
if $k$ is a field such that there exist division algebras $D$ of index
$4n$ (e.g. a number field), then  for the group $G$ with
$G(k) = \SL_m(D)$, the subgroup of reduced norm 1 of $M_m(D)$, then
$G$ is not stably rational over $k$ and $m(G) = 4n-1$.
\halm


\begin{thebibliography}{Ti1}
\bibitem[BT]-A. Borel et J. Tits, Groupes r\'eductifs,
Pub. Math. I. H. E. S. 27 (1965), 50 - 151.
\bibitem[Ch]-V. I. Chernousov, The group of congruence coefficients of a
canonical quadratic forms and stable rationality of the PSO variety,
Math. Notes 55 (1994), 413 - 416.
\bibitem[CT]-J.-L. Colliot-Th\'el\`ene, Arithm\'etique des vari\'et\'es
rationnelles et probl\`mes birationneles, Proc. I. C. M., Berkeley,
California, 1986, 641 - 653.
\bibitem[DT]-D. \v Z. \DJ okovi\'c and Nguyen Q.Thang, Conjugacy 
classes of maximal tori in simple real algebraic groups and applications,
Canadian J. Math. 46 (1994), 699 - 717.
\bibitem[MT]-Yu. I. Manin and M. Tsfasman, Rational varieties : algebra,
geometry and arithmetic, Russian Math. Surveys 41 (1986), 51 - 116.
\bibitem[M1]-A. S. Merkurjev, Generic element in $SK_1$ for simple algebras, 
K-Theory 7 (1993), 1 - 3.
\bibitem[M2]-A. S. Merkurjev, R-equivalence and rationality problem
for semisimple adjoint classical algebraic groups, 1994 (to appear).
\bibitem[O]-T. Ono, On the relative theory of Tamagawa numbers, Ann. Math. 
82 (1965), 88 - 111.
\bibitem[P]-V. P. Platonov, On the problem of rationality of spinor
varieties, Soviet Math. Dokl.
20 (1979), 1027 - 1031.
\bibitem[Se]-M. Selbach, Klassification halbeinfacher algebraischer Gruppen,
Bonn Math. Sch. 1973.
\bibitem[St]-R. Steinberg, {\it Lectures on Chevalley groups}, Yale University,
1967.
\bibitem[T1]-Nguyen Q. Thang, On the weak approximation in algebraic groups
and a theorem of Harder, Preprint 25/89, 
Karl-Weierstrass Inst. f\"ur Math., 1989.
\bibitem[T2]-Nguyen Q. Thang, On weak approximation in algebraic groups, 
Contem. Math. 131 (1992),  423 - 426.
\bibitem[T3]-Nguyen Q. Thang, On multiplicators of hermitian forms of type
$D_n$, J. Fac. Sci. Univ. Tokyo, Sec. IA, 39 (1992), 33 - 42.
\bibitem[Ti1]-J. Tits, Classification of algebraic semisimple groups, 
Proc. Symp. Pure Math. A. M. S., v. 9 (1966), 33 - 62.
\bibitem[Ti2]-J. Tits, Groupe de Whitehead des groupes alg\'ebriques
simples sur un corps (d'apr\`es V. P. Platonov et al.), S\'em. Bourbaki
29-ann\'ee, Exp. 505, 1976/1977.
\bibitem[V]-V. E. Voskresenski, Algebraic tori, Nauka, Moscow, 1977.
\bibitem[VK]-V. E. Voskresenski and A. A. Klyachko, Toroidal Fano       
varieties and root systems, Math. USSR Izvestya 24 (1985), 221 - 244.
\end{thebibliography}
\end{document}